# Complementarity Reformulations for the Optimal Design of Distributed Energy Systems with Multiphase Optimal Power Flow

*Ishanki De Mel[a], Oleksiy V. Klymenko[a], Michael Short[a*]*

[a]Department of Chemical and Process Engineering, University of Surrey, Guildford, United Kingdom GU2 7XH

*m.short@surrey.ac.uk

**Abstract**

The design of grid-connected distributed energy systems (DES) has been investigated extensively as an optimisation problem in the past, but most studies do not include nonlinear constraints associated with unbalanced alternating current (AC) power flow in distribution networks. Previous studies that do consider AC power flow use either less complex balanced formulations applicable to transmission networks, or iterative linearisations derived from local power flow solutions and prior knowledge of the design. To address these limitations, this study proposes a new algorithm for obtaining DES design decisions subject to nonlinear power flow models. The use of regularised complementarity reformulations for operational constraints that contain binary variables is proposed. This allows the use of large-scale nonlinear solvers that can find locally optimal solutions, eliminating the need for iterative linearisations and prior knowledge while improving flexibility and accuracy. DES design models with either multiphase optimal power flow (MOPF), which captures the inherent phase imbalances present in distribution networks, or balanced optimal power flow formulations (OPF) are tested using a modified version of the unbalanced IEEE EU low-voltage network. Results are also compared with a popular linear DES design framework based on direct current (DC) approximations, which proposes an infeasible operational schedule when tested with MOPF, while the fixed design alone produces the highest annualised costs. Despite the increased complexity, DES with MOPF obtains the best solution, enabling a greater integration of solar capacity where possible and reducing total annualised cost when compared to DES with OPF. The new algorithm with complementarity reformulations achieves a 19% improvement when compared with solving a bi-level model for DES with OPF, where the entire binary topology in the nonlinear model is fixed. The study therefore enables the acquisition of DES designs that can work symbiotically within distribution networks.

## 1 Introduction

With rising concern for carbon emissions associated with large-scale electricity generation and sky-rocketing energy costs, DES that incorporate small-scale renewables are becoming an increasingly attractive option for consumers. These convert traditional consumers, who rely solely on the electrical grid to meet their demands, to "prosumers", who contribute to low carbon generation by selling excess power and reducing dependency on the grid. Increased integration of DES would also allow network operators to avoid costly grid upgrades caused by rising demand [1]. Due to these benefits, most DES are connected to low-voltage distribution networks which supply end users, as it allows the consumers and network operators to establish a symbiotic relationship.

The optimal design and operation of grid-connected DES has been a focal point in literature, especially over the past decade, as properly designed DES could minimise both total costs for the consumer and carbon emissions. Optimisation models, especially Mixed-Integer Linear Programming (MILP) formulations, have been most commonly used to holistically design DES. These take technological, economical, and environmental constraints and variables into account, while being relatively easy to formulate and solve. However, detailed network constraints related to AC power flow, which are nonlinear and nonconvex, are notably absent in such models. DC approximations that are ill-suited for AC power flow constraints have been repeatedly used in MILP design frameworks, despite the advances in other facets of DES, such as the incorporation of a variety of technologies and metering regimes [2], blockchain applications for safe energy transactions [3], and hydrogen-based hybrid technologies [4].

In previous work, we identified a small-subset of studies utilising a new combined framework, labelled DES-OPF, which attempt to bridge this gap by incorporating a standalone class of models incorporating detailed AC network constraints, called OPF [5]. Most of these studies focus on DES operation or scheduling, such as [6] and [7], where the design is immutable as it is assumed to be predetermined. Very few studies have investigated the combined DES-OPF problem within a design framework, where design variables such as the capacities and locations of distributed energy resources (DERs) are influenced by power flow constraints. This is because the nonconvexity and interdisciplinary nature of the DES-OPF problem can make these models difficult to formulate and solve. Some of the existing DES-OPF studies focusing on design propose linearisations to obtain globally optimal solutions for the mixed-integer and nonlinear problem [8,9], while others use computationally expensive Mixed-Integer Nonlinear Programming (MINLP) frameworks that cannot be easily scaled [10]. The use of post-optimisation numerical checks using established power flow simulation tools has also been proposed [8]. We have extensively investigated these existing methods and compared resulting designs in a previous study [11], and propose a bi-level or two-stage method to solve the overarching MINLP problem. This work illustrated the benefits of this method, as it values practical feasibility of the design over global optimality, while being scalable for solving larger test cases. Results demonstrate that the commonly used MILP framework with DC approximations for power flow could produce the worst-performing designs, a consequence of oversimplifying network representations to achieve global optimality.

A prominent weakness in all the studies mentioned above is the use of balanced OPF formulations, which are used to analyse power flow in long-distance and high voltage transmission networks. Transmission and distribution networks employ multiphase circuits (most commonly, three-phase) for AC power flow in reality. However, the balanced OPF formulation assumes that all phases are balanced, enabling the analysis of just one representative phase. While this is a suitable assumption for high-voltage transmission networks, it is unclear whether this formulation is suitable for active distribution networks due to the inherent imbalances present in distribution networks. Multiphase Optimal Power Flow (MOPF) is a relatively new and active area of research, primarily focused on power flow in low-voltage, unbalanced, and radial distribution networks. Most studies in literature use Interior-Point methods [12,13] or Semidefinite programming relaxations [14–18] to solve the highly nonconvex problem. The few linearisation methods proposed, such as the Taylor approximation and fixed-point linearisations [19,20] rely on strictly feasible initial points that must be obtained by solving the nonconvex power flow problem using nonlinear optimisation models or simulation tools such as OpenDSS [21]. This results in an iterative procedure between nonlinear and linear models, where solving a nonlinear power flow model is a strict prerequisite. When such models interface with DES models in the combined framework where design variables remain free, these iterative linearisations could result in computationally prohibitive models, making solving such models impractical. Furthermore, using local solutions from nonlinear models or simulations to derive and solve linear models would not result in true globally optimal solutions, and only leads to an illusion of achieving global optimality.

As a consequence of this additional complexity, there is a notable absence of DES-MOPF frameworks in literature. To the best of the authors' knowledge, only one other study exists that has attempted DES-MOPF within a design framework [22]. This study employs the fixed-point linearisation mentioned above [19] in conjunction with REopt [23] for design decisions associated with rural microgrids. The limitations of the iterative linear procedure for MOPF are evident in this study, as the model cannot find feasible solutions for the overall design problem without a feasible MOPF solution precalculated from a nonlinear power flow model or simulation tool. Furthermore, a feasible MOPF initial point for the design can only be obtained if there is prior knowledge about the design and its generational capabilities at each time point. This prior knowledge of the design and associated power flow may not be available to the modeller, unless iterative pre-optimisation calculations are conducted. The cyclical iterative approach also brings to question whether the final solution of such a framework can guarantee local optimality, let alone global optimality, due to the restrictions imposed by the fixed-point linearisation. The lack of a DES-MOPF framework which can eliminate the need for prior knowledge of the design and power flow, and dependence on external power flow tools presents a significant gap in literature to be addressed.

### 1.1 Contributions of this study

The literature review above illuminates the absence of a tractable and accessible DES-MOPF framework in literature, which is essential to investigate DES design within unbalanced distribution networks. This study contributes to greater understanding of the impacts of MOPF on resulting DES designs by presenting a novel framework to solve the complex DES-MOPF problem, by extending and improving a two-stage approach we previously proposed for solving the DES-OPF problem [11]. This approach showed great promise in terms of tractability when solving a large case study. The new and improved algorithm proposed in this paper further exploits the overall structure of the combined problem by employing generic complementarity reformulations for operational constraints with binary variables. This has been inspired by the reformulation proposed by Nazir and Almassakhi [24] specific to battery operation. However, note that the reformulations employed here are more generic and reliable, as they do not require simplified battery efficiencies to predict battery operation, as done in Nazir and Almassakhi [24]. To the best of the authors' knowledge, this is the first time such reformulation techniques are employed in DES modelling. It does not rely on linear or convex approximations, which cannot guarantee the feasibility of the resulting designs unless the exactness of these approximations is verified using nonlinear power flow models or simulations. Nor does this framework rely on prior knowledge of the design or power flow solution.

As opposed to other studies in literature that use an MILP framework, this study does not focus on finding globally optimal solution. Rather, the emphasis is on finding good locally optimal solutions that minimise costs to the consumer while complying with network constraints. The framework is further detailed in Section 2, and is tested using a reduced version of the 906-node unbalanced IEEE EU LV test case [25], presented in Section 4. The DES-MOPF framework is compared with the MILP framework that employs DC approximations, and the DES-OPF framework that employs balanced power flow formulations.

## 2 Methodology

When modelling and solving optimisation problems, there are two complementary aspects to consider. The first is the optimisation framework, which contains the mathematical formulations of the physical systems subject to any simplifying assumptions. The second is the algorithm or solution method used to solve the framework. Different solution methods can be used to solve the same framework, as long as these are compatible with one another. In this paper, the algorithms for solving the DES-MOPF problem are presented first in Section 2.1. This includes the novel algorithm with complementarity reformulations for greater compatibility between the framework and solution method. The DES and MOPF formulations, along with relevant reformulations for the DES, are then introduced in Sections 2.2 and 2.3, respectively. Linking constraints that combine the DES and MOPF formulations to form the overarching DES-MOPF framework are described in Section 2.4. The objective function of the DES-MOPF problem is presented in Section 2.2.

### 2.1 Algorithms for integrating detailed power flow

Figure 1 outlines the solution methods employed in this study. The method labelled MILP is one of the most common methods in literature, where AC power flow constraints are replaced with the linear

DC power flow approximation in the formulation so that the overall framework can be solved as an MILP. When incorporating nonlinear power flow formulations (either OPF or MOPF, but especially for MOPF), Figure 1 outlines two solution methods, BL and NBL.

The previously proposed two-stage method [11], labelled BL, includes two main levels or steps. A DES formulation with the DC approximation is solved as a MILP in Level 1, followed by the addition of nonlinear power flow constraints (either OPF or MOPF) in Level 2. Both levels have a common objective, to minimise the total costs of the system. Level 1 is solved first to obtain the binary topology in which the technology locations and operational binary variables are selected. This binary topology is then fixed prior to solving the full NLP model with MOPF. This is to avoid solving an unscalable MINLP problem, as demonstrated in past work [11]. The NLP retains key continuous design decision variables such as the capacities of each technology, despite the binary topology being fixed. It also has freedom to change all continuous operational variables, thus having sufficient degrees of freedom to find a feasible solution.

The extended and improved BL method proposed in this study is labelled as NBL in Figure 1. The fixed binary topology in the second step of BL restricts the degrees of freedom, especially for key operations such as battery charging and discharging. Conventionally, binary variables have been used in DES formulations to represent complementarity problems, such as to prevent the simultaneous charging and discharging of batteries or buying and selling from the grid. While the use of binary variables to conveniently reformulate the problem to a mixed-integer and linear problem when AC power flow is not considered, these significantly increase computational complexity when combined with nonlinear AC power flow constraints. A complementarity problem involving variables $x$ and $y$ can be represented as below:

$$0 \leq x \perp y \geq 0 \tag{1}$$

and replaced by the following nonlinear complementarity constraints [26]:

$$x \cdot y = 0 \tag{2}$$

$$x \geq 0, y \geq 0 \tag{3}$$

Eq. (2) is numerically challenging to solve at the origin, and does not meet the conditions that are required by most existing large-scale NLP solvers to find feasible solutions. A commonly-used regularisation technique that reformulates Eq. (2) to improve solvability is presented below, which involves introducing a small positive parameter $\epsilon$ that is driven to zero [27]:

$$x \cdot y \leq \epsilon \tag{4}$$

The improved algorithm NBL, presented in Figure 1, replaces the eligible operational constraints with the equivalent complementarity regularisation presented in Eq. (4). The small positive parameter $\epsilon$ is driven down by a positive tuning parameter $\alpha$, until a set tolerance is met. As it is possible for the solver to terminate without finding a locally optimal solution, especially as $\epsilon$ approaches zero, another tuning parameter $\sigma$ is introduced to perturb $\epsilon$. Note that $\alpha \neq \sigma$ to ensure that $\epsilon$ does not return to a previous value. While there is no set guidance for choosing these parameters, both parameters should be greater than 1, and the choices used in this study are described further in Section 4.

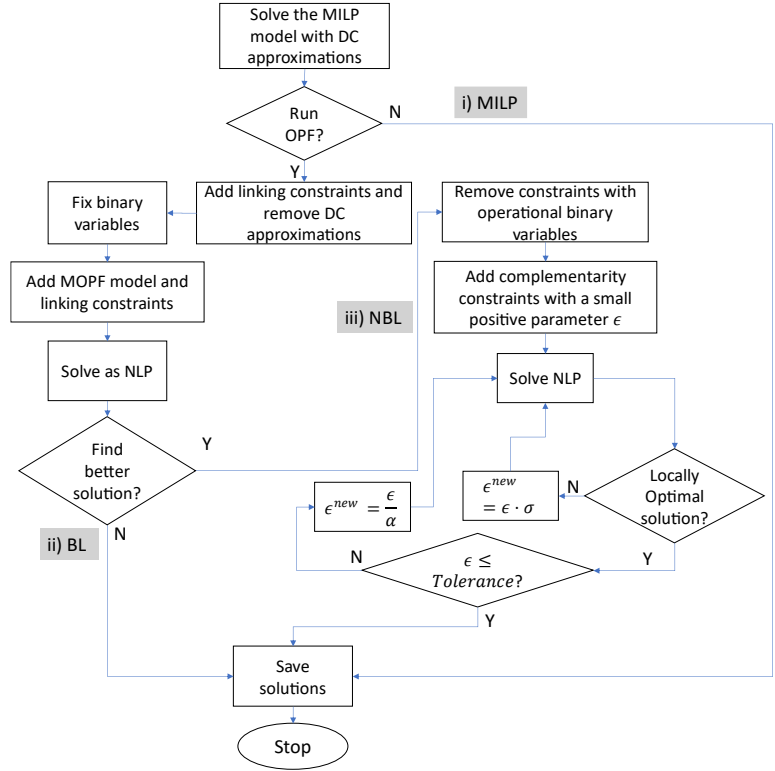

Figure 1. Overarching algorithm for solving the combined DES-MOPF problem. The problem can be solved as i) a MILP which includes linear DC power flow constraints, ii) an NLP in BL where the binary topology is fixed, and iii) the method proposed in this study, NBL, which reformulates some of the binary constraints as complementarity constraints to give the model greater degrees of freedom.

Note that some design-related binary variables may need to remain fixed in the NLP, as they may not be eligible for complementarity reformulation. These could include the discrete selection of technology types, for which reformulation may hinder the overall solution method by increasing the number of variables without significantly improving the solution.

As indicated in Figure 1, the initialisation strategy for the proposed algorithm NBL essentially involves solving BL, i.e., solving the NLP with a fixed binary topology prior to the reformulations. This ensures that a feasible initial point which respects the complementarity constraints is available to the solver when executing the NBL algorithm, one that does not involve use of external simulation tools.

The new algorithm NBL has the potential to find better solutions compared to the previous BL method, as the reformulations increase the degrees of freedom which the nonlinear constraints can influence.

## 2.2 DES formulation

The DES formulation including the DC power approximation is used, as presented in [11]. The overarching objective of both levels is to minimise the total annualised cost $TAC$ of the system across seasons $s \in S$:

$$\min \ TAC = \sum_{G \in DER} C_G^{INV} + \sum_{s \in S} \left( C_s^{grid} + \left( \sum_{G \in DER} C_{s,G}^{OM} \right) - CI_s \right) \quad (5)$$

where all the technologies considered are indicated by set $G \in \mathrm{DER}$ for notational simplicity. Costs include overall investment $C_G^{INV}$ and operation $C_{s,G}^{OM}$ of technologies, electricity purchase $C_s^{grid}$, and income $CI_s$ from both generational and export tariffs (as provided by the UK Feed-in-Tariff scheme [28]). As this formulation has the capability to optimise each season separately if necessary, the capacity of each technology at each residence or house $i \in I$ is linked across seasons using the following constraint:

$$If\ s > 1 \quad Cap^G_{i,s} = Cap^G_{i,s-1} \quad \forall i \in I \tag{6}$$

In each seasonal model, power and heating balances are used to ensure that electricity and heating demands ($E^{load}_{i,t}$ and $H^{load}_{i,t}$, respectively) are met at each time point $t \in T$, as illustrated by Eqs. (7) and (8). Note that the seasonal subscript $s$ has been excluded hereafter to simplify the notation. The electricity balance considers power from the external electrical grid $E^{grid}_{i,t}$, generation from solar photovoltaics (PVs) used to meet demand at the household $E^{PV}_{i,t}$, and discharge from lithium-ion batteries $E^{disch}_{i,t,c}$. Note that set $c \in C$ indicates types of batteries which can be included in the model. The heating balance considers heat generated from boilers $H^B_{i,t}$. The most commonly employed DERs have been chosen for this study as it focuses on the impact on existing electrical networks, and therefore these technologies are the most representative.

$$E^{load}_{i,t} = E^{grid}_{i,t} + E^{PV}_{i,t} + \sum_c E^{disch}_{i,t,c} \quad \forall i \in I, t \in T \tag{7}$$

$$H^{load}_{i,t} = H^B_{i,t} \quad \forall i \in I, t \in T \tag{8}$$

The system is further constrained to prevent the simultaneous buying and selling of electricity to and from the grid, which could lead to illicit exploitation of export tariffs [29]:

$$E^{grid}_{i,t} \leq E^{load}_{i,t} \cdot (1 - X_{i,t}) \quad \forall i \in I, t \in T \tag{9}$$

$$E^{sold}_{i,t} \leq M \cdot X_{i,t} \quad \forall i \in I, t \in T \tag{10}$$

where binary variable $X_{i,t}$ is used to linearise a complementarity problem. Note that Eq. (10) is a big-M constraint, where M is a pre-defined large value. These constraints are eligible for complementarity reformulation and regularisation employed in the NBL algorithm, where $\epsilon$ is a small positive parameter:

$$E^{grid}_{i,t} \cdot E^{sold}_{i,t} \leq \epsilon \quad \forall i \in I, t \in T \tag{11}$$

$$E^{grid}_{i,t} \leq E^{load}_{i,t} \quad \forall i \in I, t \in T \tag{12}$$

The total electricity generated from PVs can be consumed by the household, stored in batteries, and any excess can be sold to the grid. It is limited by the efficiency of the panels, available area for installation, and irradiation levels. Boiler heat generation also accounts for its efficiency and fuel consumption. Note that all key design variables, such as the capacities of solar PVs, batteries, and boilers, are formulated here as continuous variables.

Other key operational constraints include electricity storage $\mathcal{E}^{stored}_{i,t,c}$ in the batteries in kWh, given by Eqs. (13) and (14). These account for charging $E^{ch}_{i,t,c}$ and discharging $E^{disch}_{i,t,c}$ in kW, as well as their respective efficiencies, $\eta^{ch}$ and $\eta^{disch}$. The time interval is indicated by $\Delta t$.

$$if\ t = start: \mathcal{E}^{stored}_{i,t,c} = \left(E^{ch}_{i,t,c} \cdot \eta^{ch} \cdot \Delta t\right) - \frac{E^{disch}_{i,t,c} \cdot \Delta t}{\eta^{disch}} \quad \forall i \in I, t \in T, c \in C \tag{13}$$

$$else: \mathcal{E}^{stored}_{i,t,c} = \mathcal{E}^{stored}_{i,t-1,c} + \left(E^{ch}_{i,t,c} * \eta^{ch} * \Delta t\right) - \frac{E^{disch}_{i,t,c} \cdot \Delta t}{\eta^{disch}} \quad \forall i \in I, t \in T, c \in C \tag{14}$$

Note that the storage term is absent in Eq. (13) as it is assumed that the initial storage level at the first timestep, $t = start$, is zero. The energy stored at the batteries is further constrained by the State of Charge $SoC^{max}_c$ and maximum depth of discharge $DoD^{max}_c$:

$$\mathcal{E}_{i,t,c}^{stored} \leq Cap_{i,c}^{batt} \cdot \text{SoC}_{c}^{\max} \quad \forall i \in \text{I}, t \in \text{T}, c \in \text{C} \tag{15}$$

$$\mathcal{E}_{i,t,c}^{stored} \geq Cap_{i,c}^{batt} \cdot (1 - \text{DoD}_{c}^{\max}) \quad \forall i \in \text{I}, t \in \text{T}, c \in \text{C} \tag{16}$$

The batteries are also allowed to charge from power supplied by the external electrical grid. However, they are not allowed to charge and discharge at the same time, which is governed by the following big-M constraints:

$$E_{i,t,c}^{ch} \leq \text{M} \cdot Q_{i,t,c} \tag{17}$$

$$E_{i,t,c}^{disch} \leq \text{M} \cdot (1 - Q_{i,t,c}) \tag{18}$$

where $Q_{i,t,c}$ is a binary variable. These constraints are also eligible for complementarity reformulation and regularisation in the NBL algorithm, eliminating the binary variable, as shown below:

$$E_{i,t,c}^{ch} \cdot E_{i,t,c}^{disch} \leq \epsilon \tag{19}$$

Although the type of battery is restricted to Lithium-ion in this study, it is possible to present more than one battery type as options in the model, for e.g., lithium-ion, lead-acid, etc. These battery types can be selected using the binary variables $W_{i,c}$. However, only one battery type should be installed at each dwelling. These conditions are represented by the constraints below:

$$\sum_{c} W_{i,c} \leq 1 \tag{20}$$

$$Cap_{i,c}^{batt} \leq \text{M} \cdot W_{i,c} \tag{21}$$

As there is no straightforward complementarity reformulation applicable to these constraints, the binary variable $W_{i,c}$ would be determined by the MILP and remain fixed when power flow constraints are introduced in the NLP.

The DC power flow approximation is presented below, which considers only active power injections $P_{n,t}$ at node $n \in \text{N}$ of branch $(n,m) \in \text{L}$:

$$P_{n,t} \approx \sum_{m=1}^{N} \text{B}_{nm}(\theta_{n,t} - \theta_{m,t}) \quad \forall n \in \text{N}, t \in \text{T} \tag{22}$$

where the parameter $\text{B}_{nm}$ represents susceptance, while $\theta_{n,t}$ represents voltage angle. Note that susceptance is obtained from the single-phase or balanced version of admittance calculations, which is further explained in [11]. The term for voltage magnitude, which is found in the nonlinear AC power flow formulations, is not included in this approximation as it assumes nominal voltage at each node. Active power losses, which are captured by the differences between voltage angles and magnitudes, and reactive power calculations are all ignored. The assumptions resulting in the DC approximation hold for transmission networks as they typically have a low $R/X$ ratio, but are unrealistic for multiphase AC power flow in low-voltage distribution networks which typically have a high reactance-to-resistance ($R/X$) ratio.

### 2.3 MOPF formulation

There are several different methods for formulating (and solving) the MOPF problem, such as iterative Z-bus [20], current injection [12,17], and branch flow [30]. The phase-decoupled bus injection model is utilised here as this formulation is the most compatible MOPF formulation when interfacing with DES. It enables all buses to be modelled as constant active power-reactive power ($PQ$) buses, which in turn allows power injections at load/generational nodes that are linked to the DES design formulation, described further in Section 2.4. This formulation has been used in multiphase state estimation models, such as [31]. The per-unit system is used for tractability, where line-to-line voltages $V_{LL}^{base}$ are chosen as bases for both the source (primary) side and load (secondary) side of the transformer, along with an apparent power base $S^{base}$ common to both sides [32].

The main power flow equations for active power $P_{n,t}^{\phi}$ and reactive power $Q_{n,t}^{\phi}$ at node $n \in N$, timepoint $t \in T$, and phase $\phi \in \Phi$ are presented below:

$$P_{n,t}^{\phi} = V_{n,t}^{\phi} \sum_{m \in N} \sum_{\varphi \in \Phi} V_{m,t}^{\varphi} \left( g_{mn}^{\phi\varphi} \cos\left(\theta_{n,t}^{\phi} - \theta_{m,t}^{\varphi}\right) + b_{mn}^{\phi\varphi} \sin\left(\theta_{n,t}^{\phi} - \theta_{m,t}^{\varphi}\right) \right) \quad \forall n \in N, t \in T, \phi \in \Phi \tag{23}$$

$$Q_{n,t}^{\phi} = V_{n,t}^{\phi} \sum_{m \in N} \sum_{\varphi \in \Phi} V_{m,t}^{\varphi} \left( g_{mn}^{\phi\varphi} \sin\left(\theta_{n,t}^{\phi} - \theta_{m,t}^{\varphi}\right) - b_{mn}^{\phi\varphi} \cos\left(\theta_{n,t}^{\phi} - \theta_{m,t}^{\varphi}\right) \right) \quad \forall n \in N, t \in T, \phi \in \Phi \tag{24}$$

where $V_{n,t}^{\phi}$ represents the voltage magnitude, $\theta_{n,t}^{\phi}$ is the voltage angle, as the voltage is represented in polar form. For a three-phase network, the set of phases would be defined as $\phi \in \Phi = \{a, b, c\}$. Note that $m \in N$ also represents buses or nodes, while $\varphi \in A$ also represents phases. The terms $g_{mn}^{\phi\varphi}$ and $b_{mn}^{\phi\varphi}$ are Conductance G and Susceptance B, respectively, which are obtained from the three-phase admittance matrix, Y, as shown in Eq. (25) below [33]. The calculation of Y, as outlined in Eqs. (25)-(28), has been included in the overall optimisation framework, such that these parameters do not need to be obtained via external power flow simulation tools. Note that conductance and susceptance are treated as parameters in this optimisation model as the network configuration is immutable (i.e., there are no reconfiguration capabilities), representing the existing underlying distribution network to which the residential DES are connected.

$$Y = \begin{bmatrix} Y_{11}^{abc} & \cdots & Y_{1m}^{abc} \\ \vdots & \ddots & \vdots \\ Y_{n1}^{abc} & \cdots & Y_{nm}^{abc} \end{bmatrix} = G + jB \tag{25}$$

For a three-phase network, each element of Y represents a $3 \times 3$ matrix. These submatrices can be obtained by applying the Approximate Line Model [32] using parameters for resistance R and reactance X for each line in the network. Then, similar to balanced OPF, the diagonal and non-diagonal elements of Y can be obtained from the equations below (adapted from [34] for a three-phase network):

$$Y_{nn}^{abc} = \sum_{m:(n,m) \in L} Y_{nm}^{abc} + \sum_{m:(m,n) \in L} Y_{mn}^{abc} \tag{26}$$

$$Y_{nm}^{abc} = - \sum_{m:(n,m) \in L} Y_{nm}^{abc} - \sum_{m:(m,n) \in L} Y_{mn}^{abc} \quad n \neq m \tag{27}$$

Note that line shunt admittances have been excluded in the above equations as they are typically ignored for short overhead and underground distribution cables. The three-phase admittance matrix $Y^T$ for the transformer, taking primary ($p$) and secondary ($s$) taps into account, are also included in the above calculations [33]:

$$Y^T = \begin{bmatrix} Y_{pp}^{abc} & Y_{ps}^{abc} \\ Y_{sp}^{abc} & Y_{ss}^{abc} \end{bmatrix} \tag{28}$$

Each element in the $Y^T$ matrix is also a $3 \times 3$ matrix, whose calculations are based on the original transformer submatrices proposed by Chen et al. [35]. These calculations account for the copper and core losses, winding connections, phase shifting, and off-nominal tapping. Note that a small shunt admittance has been included in the matrix calculations to avoid singularities [36].

In grid-connected DES and microgrids, the bus that acts as the point of common coupling, i.e. the bus at which the distribution network connects to the wider external network, is treated as a reference bus for voltage. This reference bus is indicated by subscript $ref$, and is often also used as the slack bus (indicated by slack). A slack bus allows unlimited active and reactive power injections, representing the power injections from the external grid. As the source or reference is typically at primary side of the transformer, therefore having the slack bus here would present the power injections at the primary side. The voltage magnitude at the reference bus is fixed to the per-unit voltage prescribed at the source:

$$V_{ref,t}^{\phi} = V_{\text{source}} \quad \forall t \in T, \phi \in \Phi \tag{29}$$

The voltage angles at the reference bus are also fixed to those of the source, as demonstrated for a three-phase network below:

$$\theta_{ref,t}^{a} = 0 \quad t \in T$$
$$\theta_{ref,t}^{b} = -120° \quad t \in T \tag{30}$$
$$\theta_{ref,t}^{b} = 120° \quad t \in T$$

The voltage magnitude at each node and timepoint must remain within pre-specified network bounds:

$$\frac{V^{UB}}{V_{\text{base}}} \geq V_{n,t}^{\phi} \geq \frac{V^{LB}}{V_{\text{base}}} \quad \forall n \in N, t \in T, \phi \in \Phi \tag{31}$$

The constraints linking the power injections for generational and load nodes, defined by set $n \in I \subseteq N$, are described further in Section 2.4 below, which completes the MOPF formulation. Nodes or buses that act as connectors do not have any active or reactive power injections. These buses are defined as $n \notin \{I, \text{slack}\}$ and are governed by the following constraints:

$$P_{n,t}^{\phi} = 0 \quad \forall n \notin \{I, \text{slack}\}, t \in T, \phi \in \Phi \tag{32}$$

$$Q_{n,t}^{\phi} = 0 \quad \forall n \notin \{I, \text{slack}\}, t \in T, \phi \in \Phi \tag{33}$$

## 2.4 Linking DES and MOPF

The net active and reactive power injections at the load nodes, defined by set $n \in I \subseteq N$ (where index $i \in I$ is used in the DES formulation), are the main interfacing points between the DES and MOPF. The net injections at connecting nodes $n \notin \{I, \text{slack}\}$, which do not consume nor generate power, must equal zero. Loads in unbalanced multiphase networks are not equally distributed between three phases, constraints linking the DES and MOPF need to consider the phase to which the load is connected. Typically, residential loads are connected to one out of the three phases, especially in radial distribution networks. As this information is unavailable to the DES formulation, a parametric indicator $\Psi_n^{\phi} \in \{0,1\}, \phi \in \Phi$ is used to indicate the phase to which each house $i$ is connected, as shown in the linking constraints below.

$$P_{n,t}^{\phi} = \Psi_n^{\phi}\left(P_{n,t}^{Gen} - P_{n,t}^{Load}\right) = \Psi_n^{\phi}\left(\frac{E_{i=n,t}^{PV,sold} - E_{i=n,t}^{grid} - \sum_c E_{i=n,t,c}^{grid,charge}}{S^{\text{base}}}\right) \forall n \in I, t \in T, \phi \in \Phi \tag{34}$$

$$Q_{n,t}^{\phi} = \Psi_n^{\phi}\left(Q_{n,t}^{Gen} - Q_{n,t}^{Load}\right) = \Psi_n^{\phi}\left(\frac{Q_{i=n,t}^{Gen} - Q_{i=n,t}^{Load}}{S^{\text{base}}}\right) \forall n \in I, t \in T, \phi \in \Phi \tag{35}$$

Note that $E_{i=n,t}^{PV,sold}$ is the excess energy sold to the grid, while $Q_{i=n,t}^{Gen}$ and $Q_{i=n,t}^{Load}$ are reactive power generation and load, respectively. $S^{base}$ is the base apparent power used to convert the powers to the per unit system (see Section 2.3). In this formulation, reactive power generating technologies are not considered, and a constant power factor is assumed to calculate the reactive power loads.

## 3   Case Study

To enable extensive testing and comparison, a modified version of the 906-node IEEE EU LV test case [25] is utilised, with 22 loads and 389 buses/nodes in total. It includes a three-phase Delta-Wye ($\Delta - Y$) transformer with a -30° phase shift. All the nodes are connected to one out of the three phases, as illustrated in Figure 2. The network is considered unbalanced because the load capacities in each of the phases are unequal. The slack bus and reference bus are both the primary side of the transformer, labelled Bus 0.

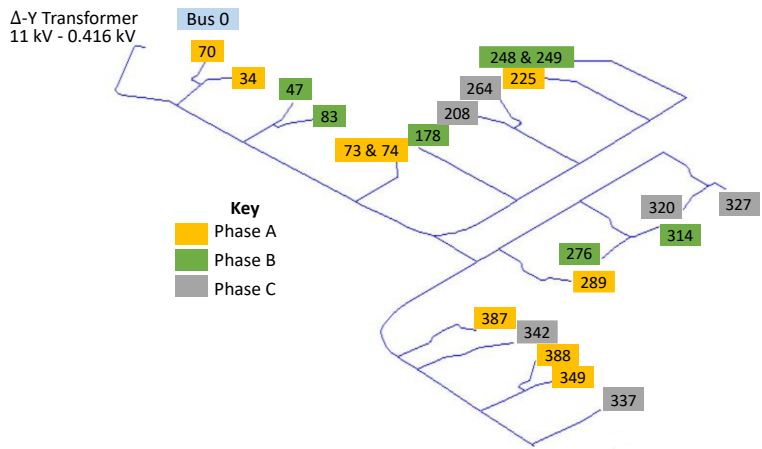

Figure 2. The modified version of the IEEE EU LV test case with 22 loads.

The test case provides an electricity demand profile for one day in winter in 1-minute intervals. This is converted to four seasonal profiles with hourly intervals, by obtaining hourly averages from the original single-day profile and scaling the overall profile for each season based on trends observed in other residential datasets [37]. As heating profiles are not provided in the original test case, these are generated based on a residential profile from Morvaj *et al.* [8]. Seasonal averaged weather profiles have been obtained from data for Charlwood, Surrey, U.K., available from the Centre for Environmental Data Analysis (CEDA) [38]. Figure 3 portrays the seasonal electricity and heating profiles for Bus 34, as an example.

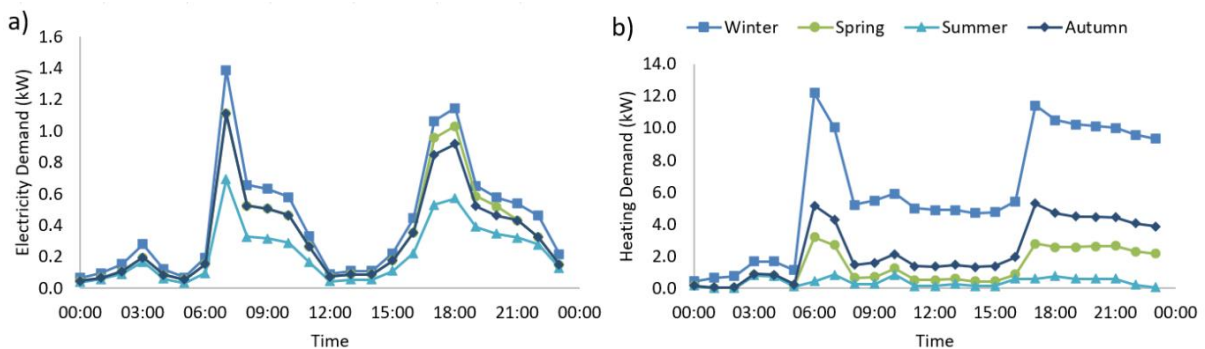

Figure 3. The seasonal profiles of a) electricity demand and b) heating demand of Bus 34.

The Economy 7 tariff [39] has been used, which has a higher electricity price of £0.18/kWh during daytime and a night price of £0.08/kWh [40]. A gas price of £0.02514/kWh is used, obtained from private correspondence. It is assumed that all the load buses are eligible for DER installation, with 35 m$^2$ available roof area for PV installation and 0.5 m$^3$ space available for battery installation. The DES total lifetime is 20 years, in line with the fixed tariffs provided by the FIT scheme [28], with an assumed interest rate of 7.5%. Under the FIT scheme, an export tariff of £0.0503/kWh is used [41], while a generational tariff of £0.1/kWh is assumed (which is higher than nominal values provided) to encourage PV installation. Table 1 presents key parameters associated with the technologies used. Table 2 shows the average and peak demands of the loads in winter, which is the season with the highest demand. All the profiles, data, and source code are provided in the following Github repository: https://github.com/Ishanki/DES-Multiphase-OPF.

Table 1. Technology parameters used in the case study.

| Technology | Parameter | Value | Unit | Reference |
|---|---|---|---|---|
| PVs | PV Investment cost | 450 | £/panel | Average from [42] |
|  | PV Efficiency | 0.18 | - | [43] |
|  | PV Fixed operational cost | 12.5 | £/kW-yr | [43] |
|  | Panel area | 1.75 | m$^2$ | Approximation from [42] |
|  | Panel capacity | 0.25 | kW | Average from [42] |
| Boilers | Boiler investment cost | 120 | £/kW | [44] |
|  | Boiler efficiency | 0.94 | - | [45] |
| Batteries | Volumetric energy density | 20 | kWh/m$^3$ | [43] |
|  | Max DoD | 0.85 | - | [43] |
|  | Max SoC | 0.9 | - | [43] |
|  | Investment cost | 270 | £/kWh | [43] |
|  | Operational cost | 11 | £/kWh-yr | [43] |
|  | Round trip efficiency (RTE) | 0.89 | - | [43] |
|  | Charge efficiency | 0.94 | - | Calculated using RTE |
|  | Discharge efficiency | 0.91 | - | Calculated using RTE |

Table 2. Average and peak demands of load buses in the network in winter (the season with the highest demand).

| Load | Phase | Average Demand (kWh) | Peak Demand (kWh) |
|---|---|---|---|
| 34 | A | 0.43 | 1.39 |
| 47 | B | 0.50 | 2.70 |
| 70 | A | 0.28 | 1.20 |
| 73 | A | 0.41 | 1.58 |
| 74 | A | 0.29 | 1.08 |
| 83 | B | 0.22 | 0.89 |
| 178 | B | 0.42 | 1.88 |
| 208 | C | 0.65 | 1.83 |
| 225 | A | 0.81 | 2.86 |
| 248 | B | 0.47 | 2.44 |
| 249 | B | 0.24 | 1.08 |
| 264 | C | 0.23 | 0.51 |
| 276 | B | 0.43 | 2.17 |
| 289 | A | 0.22 | 0.63 |
| 314 | B | 0.26 | 1.01 |
| 320 | C | 0.44 | 1.28 |
| 327 | C | 0.22 | 0.72 |
| 337 | C | 0.42 | 1.85 |

| 342 | C | 0.46 | 1.52 |
| 349 | A | 0.54 | 1.97 |
| 387 | A | 0.22 | 1.23 |
| 388 | A | 0.20 | 0.60 |

## 4 Results and Discussion

Table 3 summarises the different solution methods used for testing and comparison, and their compatibility with nonlinear OPF and MOPF formulations for the combined framework. The MILP uses the DES formulation with DC approximations. Solutions from the MILP are further tested using MILP Checks, using either BL or NBL, which include MOPF formulations. The values of the continuous design variables for capacities of technologies obtained from the MILP are fixed in MILP Checks, alongside the binary topology. Note that the operational continuous variables are left as free variables, firstly because further checks confirm that fixing these variables results in infeasible solutions with respect to OPF and MOPF, and secondly, to fix the proposed designs and test their feasibility while operation is allowed to vary, as is the case in reality.

Table 3. The solution methods tested in this study.

| Solution Method | Design continuous variables | Design binary variables | Operational continuous variables | Operational binary variables | Nonlinear OPF/MOPF |
|---|---|---|---|---|---|
| MILP | Free | Free | Free | Free | N |
| BL | Free | Fixed | Free | Fixed | Y |
| NBL | Free | Fixed | Free | Free | Y |
| MILP Check (using BL) | Fixed | Fixed | Free | Fixed | Y |
| MILP Check (using NBL) | Fixed | Fixed | Free | Free | Y |

The solvers employed to test the methods proposed in Section 2.1 are CPLEX [46] for MILP steps, and CONOPT [47] for NLP steps. All tests were conducted on an Intel® Core™ i7-10510U CPU at 1.80GHz - 2.30 GHz. The MILP with DC approximations contains 134,519 constraints, 68,461 continuous variables, and 4,312 binary variables. The DES-OPF model contains 254,911 constraints and 166,573 continuous variables post-reformulation of operational binary variables. The DES-MOPF framework contains 476,191 constraints and 387,757 continuous variables, post-reformulation.

As mentioned in Section 2.1, the NBL algorithm's initial value for $\epsilon$ was taken as 0.1. The parameter $\alpha$ was chosen to be 10, such that sufficient progress is achieved with each iteration. The perturbation parameter $\sigma$ was chosen to be 15, such that if the reduction of $\epsilon$ returned a solution other than locally optimal, $\epsilon$ was increased to a value greater than the product of $\alpha$ and its previous value. The tolerance set for the complementarity regularisation parameter $\epsilon$ was $1 \times 10^{-7}$. For this test case, the magnitudes of the variables entering the complementarity constraints predominantly tend to be in the order of $0.01 - 0.1$, and the tolerance ensures complementarity by returning extremely small values (or zero for most cases) for one variable.

### 4.1 Results from the proposed algorithm (NBL)

The results for the frameworks/formulations solved using the NBL method are presented in Table 4. The MILP results (using the DC approximation) are presented in the same table for comparison. The MILP Check assessing the feasibility of the MILP design uses the MOPF formulation, as it is most representative of the test case network. The DES-OPF framework includes nonlinear balanced optimal power flow constraints, while the DES-MOPF includes multiphase optimal power flow constraints.

Predictably, the MILP obtains the best objective value. However, the MILP Check with the fixed design affirms that this optimistic objective value is unattainable with respect to multiphase network constraints, due to the reduction in income from energy generation and export as a result of renewable energy curtailment. This is because the DC approximations cannot accurately represent the power flows of AC distribution networks, and the MILP Check reveals that the design leads to voltage violations during baseline operation. The DES-MOPF formulation presents the lowest total annualised cost while respecting power flow constraints, outperforming the DES-OPF formulation. This is further investigated in Figure 4, where the capacities of the PVs and batteries at each bus are compared.

Table 4. Summary of results and costs for the combined DES-power flow formulations, solved using the NBL algorithm.

|  | **MILP** | **MILP Check** | **DES-OPF** | **DES-MOPF** |
| --- | --- | --- | --- | --- |
| Objective (£) | 21,297 | 23,103 | 23,484 | 22,072 |
| % Difference (w. MILP) | - | 8 | 10 | 4 |
| Grid Cost (£) | 4,962 | 4,962 | 5,003 | 4,265 |
| PV Investment (£) | 19,422 | 19,422 | 8,477 | 15,424 |
| PV Operation (£) | 1,375 | 1,375 | 600 | 1,092 |
| Boiler Investment (£) | 3,359 | 3,359 | 3,359 | 3,359 |
| Boiler Operation (£) | 13,804 | 13,804 | 13,804 | 13,804 |
| Battery Investment (£) | 64 | 64 | 624 | 565 |
| Battery Operation (£) | 27 | 27 | 259 | 235 |
| Export Income (£) | 6,368 | 5,763 | 2,013 | 4,543 |
| Generation Income (£) | 15,348 | 14,147 | 6,628 | 12,127 |
| Time taken (s) | 82 | 3,434 | 1,151 | 3,732 |
| Solver termination | Optimal | Locally Optimal | Locally Optimal | Locally Optimal |

It is evident from Figure 4 that the MILP chooses the highest PV capacities and lowest battery capacities, while DES-OPF takes a more conservative approach by choosing the highest battery capacities and lowest PV capacities. This is because the approximations in the MILP do not detect voltage violations at all, while the DES-OPF penalises most generating nodes due to the assumption that all phases are balanced and therefore voltage violations are much more prominent. The DES-MOPF takes a less conservative approach by increasing PV capacities, with battery capacities remaining close to those predicted by the DES-OPF. This is because the DES-MOPF formulation assesses the imbalances of each phase, and therefore can more accurately determine the voltage violations compared to both MILP and DES-OPF.

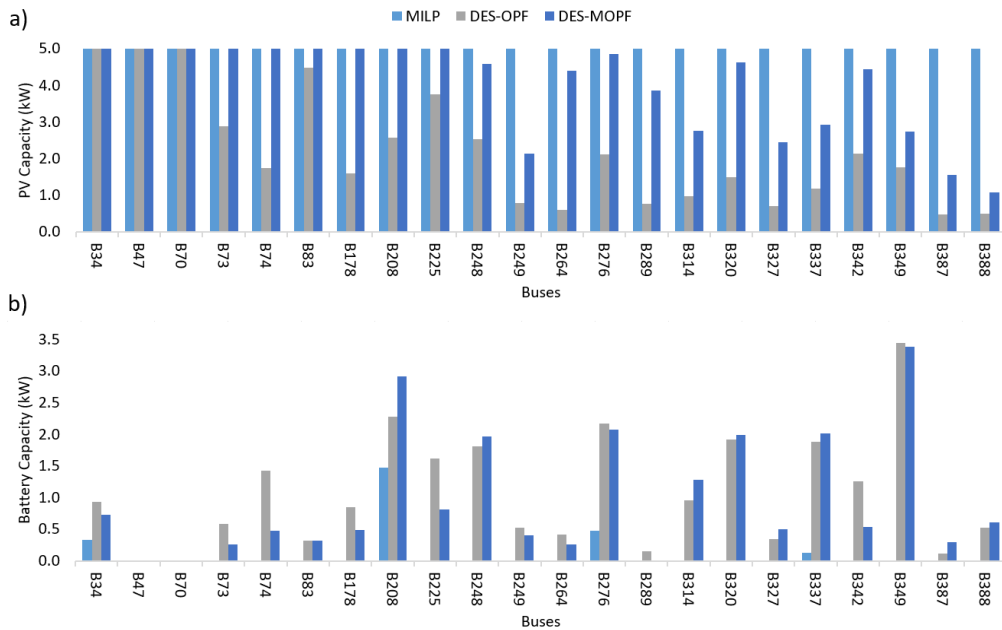

Figure 4. The resulting capacities of a) PVs and b) batteries across the three formulations: MILP, DES-OPF, and DES-MOPF.

It also appears from Figure 4 that the PV capacity of buses located towards the latter part of the network gets reduced by DES-OPF and DES-MOPF when compared to the MILP result. This could potentially be due to network violations becoming more pronounced towards the latter part of the network, especially for phases A and C, as evident in Figure 5. This figure also demonstrates how the voltage at these buses increases at the times when renewable energy export occurs. The effects on battery capacity are therefore potentially impacted by the location of the load bus on the network, as well the peak demand of the loads. The buses with the highest battery capacity installed have relatively high peak demands, though not the highest average or peak demands. Overall, the impacts on both PV and battery capacities are not straightforward to predict due to the influence of voltage unbalances and violations across the phases. These results suggest that the most suitable framework for designing DES connected to distribution networks is the DES-MOPF framework, as it neither overestimates nor underestimates the costs and violations. They further demonstrate that globally optimal solutions obtained from conventional MILP frameworks are not always executable in practice, as the oversimplified formulation is unable to detect any network violations and grossly underestimates the battery capacities required for smoother baseline operation. Unfortunately, despite the DES-MOPF solution having a more favourable total annualised cost, the computational expense does increase significantly, especially as the DES-MOPF formulation takes nearly three times as long as the DES-OPF solution.

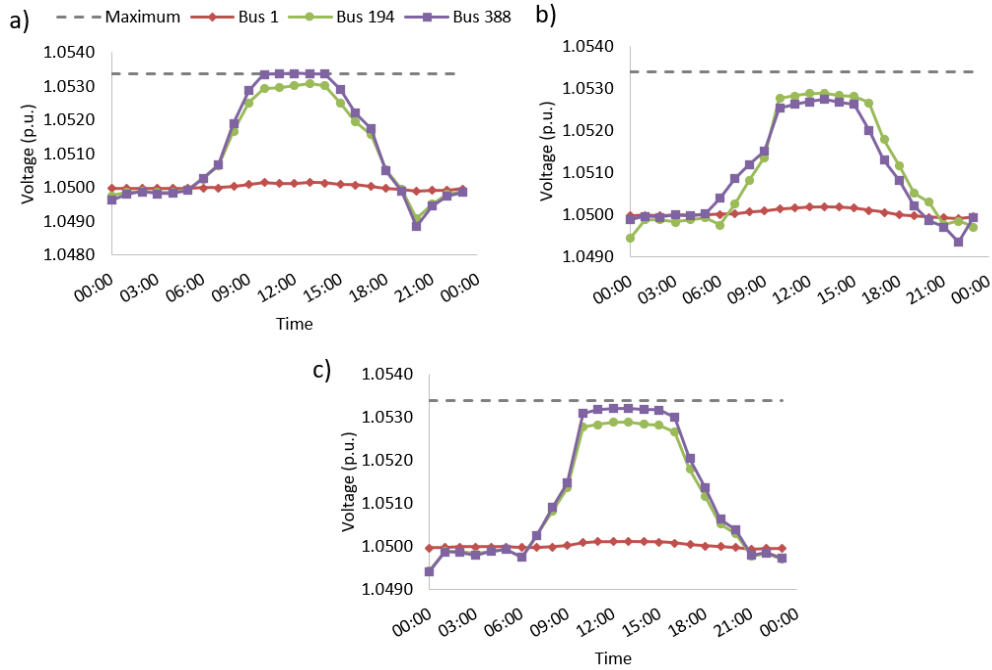

Figure 5. Voltages for Buses 1, 194, and 388 in a) phase A, b) phase B, and c) phase C, as predicted by the DES-MOPF formulation. Note that the "maximum" represents the voltage upper bound permitted by the network.

Note that the formulation of key design decisions (such as equipment capacities) as continuous variables, as commonly done in DES design models, allows the nonlinear part of the proposed algorithm to influence them and obtain better solutions overall. This would not be possible if these decisions are made discrete, and would pose a limitation to using the algorithm and combined framework, as the power flow constraints would now have limited influence on the design. The impact from making such discrete decisions on the results produced by the framework should be investigated and addressed in future work.

### 4.2 Comparison of methods

The performance of the algorithms NBL and BL are compared by calculating the percentage differences between NBL and BL objective values for each framework, DES-OPF, DES-MOPF, and the MILP Check using MOPF. These results are summarised in Table 5. The objectives for MILP Check have a negligible percentage difference between the two algorithms, potentially because the design is fixed and the degrees of freedom in this model is severely restricted as a consequence. This, however, further confirms that the original MILP with DC approximations underestimates the total annualised costs for the DES. Interestingly, the objective for DES-OPF using BL is 19% higher in comparison with NBL for the DES-OPF formulation, suggesting that the proposed algorithm is indeed capable of finding better solutions, especially due to the increased degrees of freedom provided by the reformulations. The solution for DES-MOPF produced from the NBL algorithm is also better compared to the BL algorithm, though not significantly so. The DES-MOPF formulation inherently contains higher degrees of freedom when compared to the DES-OPF formulation, due to modelling of separate phases as opposed to a single balanced phase. Therefore, it is possible that the additional degrees of freedom provided by the binary reformulations do not have as much an impact when solving the DES-MOPF formulation. These results suggest that the NBL algorithm is essential when solving DES-OPF formulations to find better feasible solutions, but the BL algorithm may be sufficient when solving DES-MOPF formulations or testing the MILP designs with respect to multiphase power flow constraints. The improvements in objective value offered by the NBL algorithm come at a higher computational expense, where it takes approximately twice as long as the BL algorithm, and therefore using BL for the DES-MOPF framework could potentially produce a good local solution in less time.

Table 5. Comparison of algorithms BL and NBL to assess the solution quality of nonlinear model formulations.

|  | MILP Check | DES-OPF | DES-MOPF |
|---|---|---|---|
| NBL Objective (£) | 23103 | 23484 | 22072 |
| BL Objective (£) | 23107 | 27988 | 22669 |
| % Difference | 0.02 | 19.18 | 2.70 |

### 4.3 Verification of the MOPF results

The MOPF formulation presented in Section 2.3 was verified by simulating the original IEEE EU LV Test feeder [25] without any DER installations, with only the data provided from the test case. The resulting voltages were then compared to the GridLab-D [48] solutions across all nodes and phases for the same test case and for three specified time intervals or snapshots: Off-peak (minute 1), On Peak (minute 566), and End (minute 1440). No voltage upper bound was considered in the MOPF simulation, as the GridLab-D simulation does not impose constraints on voltage as normally done in optimal power flow problems. The maximum and minimum percentage differences observed are reported in Table 6. The On Peak simulation reported a -4% error, which is the timestamp with the highest demand. This could potentially be a result of the different solution methods employed by GridLab-D and the MOPF simulation. As the MOPF simulation has negligible percentage errors with respect to the other two timestamps, and a comparable percentage error for the timestamp which causes the greatest strain on the network, the MOPF is deemed sufficiently accurate for use within the DES-MOPF framework.

Table 6. Comparison of MOPF simulation results with GridLab-D simulations for the IEEE EU LV test case.

|  | Off Peak | On Peak | End |
|---|---|---|---|
| Min (%) | -0.073 | -3.930 | -0.294 |
| Max (%) | 0.023 | 0.673 | 0.082 |

## 5 Conclusions

The study presents a combined framework for optimising the designs of DES connected to low-voltage distribution networks, while simultaneously considering nonlinear multiphase optimal power flow. Mixed-integer and linear economic, operational, and design constraints for the DES are considered, alongside detailed nonlinear power flow constraints representing the distribution network through which the DES imports and exports power. A new algorithm for obtaining reliable local solutions for the combined framework is proposed. It reformulates DES constraints with operational binary variables to nonlinear complementarity constraints, and uses a regularisation technique to find locally optimal solutions. This method is compared to a two-stage or bi-level algorithm, where the entire binary topology is first determined by solving an MILP and subsequently fixed, prior to solving the large-scale nonlinear model. Both algorithms do not rely on external power flow simulations, posing a significant advantage by allowing the power flow constraints to influence DES design decisions. The need for iterative fixed-point linearisations and prior knowledge on the design and power flows have also been eliminated.

Results of the DES model with multiphase optimal power flow are compared to alternative power flow formulations that have been frequently employed in DES literature, such as the MILP model with the DC approximations, and the nonlinear balanced power flow formulation that evaluates power flow over a single phase. These results and tests confirm that the MILP with the DC approximations underestimates the total costs due to loss of export income from renewable energy curtailment, despite the highest integration of renewable energy capacity. Furthermore, the operational decisions proposed by the MILP are infeasible when tested with multiphase power flow constraints. On the other hand, the DES model with balanced power flow overestimates the network violations and significantly reduces renewable energy capacity, resulting in an overestimation of annualised costs. The DES model with multiphase power flow constraints finds a less conservative locally optimal solution which enables greater renewable energy integration, while mitigating the power flow violations at each phase. This signifies the importance of using multiphase power flow when designing DES that are connected to low voltage distribution networks. When comparing the algorithms for solving the combined framework, the proposed algorithm produces the best solutions for the models with nonlinear power flow, with a significant improvement observed in the objective value of the DES model with balanced power flow. These results suggest that the proposed algorithm could be a valuable tool for solving the combined problem, especially when designing DES that are connected to medium-to-high voltage transmission networks where the balanced power flow assumption holds. Computational expense remains a limitation; however, the algorithm has proven its tractability for the 389-bus test case considered. Future work involves improving the algorithm to capture discrete design decisions for technologies, such as discrete capacity choices that would impact the performance of the algorithm, and assessing the tractability of the proposed algorithm for a more holistic DES.

## Acknowledgement

We thank Dr. Kyri Baker for recommending the use of the per-unit system in the three-phase formulation, which improved the overall tractability of the model.